\documentclass[10pt]{amsart}			

\usepackage{amssymb, amsmath, amsthm, mathrsfs}

\usepackage[UglyObsolete]{diagrams}

\copyrightinfo{~2012}{Bram Mesland}
\author{Bram Mesland}

\title{Spectral Triples and $KK$-Theory: A Survey}

\date\today
\address{Mathematics Institute \\
Zeeman Building\\
University of Warwick\\
Coventry CV4 7AL}
\email{brammesland@gmail.com}
\theoremstyle{plain}

\newtheorem{theorem}{Theorem}[subsection]

\newtheorem{proposition}[theorem]{Proposition}

\newtheorem{lemma}[theorem]{Lemma}

\theoremstyle{definition}

\newtheorem{definition}[theorem]{Definition}

\newtheorem{example}[theorem]{Example}
\newtheorem{remark}[theorem]{Remark}

\DeclareFontFamily{OT1}{pzc}{}
\DeclareFontShape{OT1}{pzc}{m}{it}{<-> s * [1.20] pzcmi7t}{}
\DeclareMathAlphabet{\mathpzc}{OT1}{pzc}{m}{it}

\newcommand{\Z}{\mathbb{Z}}

\newcommand{\K}{\mathbb{K}}
\newcommand{\C}{\mathbb{C}}

\newcommand{\Aut}{\textnormal{Aut}\,}

\newcommand{\Dom}{\mathfrak{Dom}\,}

\newcommand{\tildeotimes}{\tilde{\otimes}}

\newcommand{\minotimes}{\overline{\otimes}}

\newcommand{\lrh}{\leftrightharpoons}

\newcommand{\End}{\textnormal{End}}
\newcommand{\Endst}{\textnormal{End}^{*}}

\newcommand{\Sob}{\textnormal{Sob}}
\begin{document}
\maketitle
\section*{Introduction}
This survey covers the material in the forthcoming paper \cite{Mes}, which deals with the construction of a category of spectral triples that is compatible with the Kasparov product in $KK$-theory (\cite{Kas}). These notes serve as an intuitive guide to the results described there, avoiding the necessary technical proofs. We will also add some background and a broader perspective on noncommutative geometry. The theory described shows that, by introducing a notion of smoothness on unbounded $KK$-cycles, the Kasparov product of such cycles can be defined directly, by an algebraic formula. This allows one to view such cycles as morphisms in a category whose objects are spectral triples.\newline

We will consider all $C^{*}$-algebras to be equipped with a spectral triple that is sufficiently smooth. A smooth $KK$-cycle for a pair of such $C^{*}$-algebras $(A,B)$ is a triple $(\mathpzc{E},S,\nabla)$, where the pair $(\mathpzc{E},S)$ is a $KK$-cycle in the sense of Baaj-Julg \cite{BJ}, satisfying some smoothness conditions compatible with the given spectral triples, and $\nabla$ is a connection on the module, compatible with the operator $S$ and the smooth structure on $\mathpzc{E}$. Composition of such triples is defined by
 \[(\mathpzc{E},S,\nabla)\circ(\mathpzc{F},T,\nabla '):=(\mathpzc{E}\tildeotimes_{B}\mathpzc{F},S\otimes 1+ 1\otimes_{\nabla} T, 1\otimes_{\nabla}\nabla '),\]
and preserves all smoothness conditions. Moreover, it represents the Kasparov product of the $KK$-cycles $(\mathpzc{E},S)$ and $(\mathpzc{F},T)$.\newline

In particular this allows one to compute such products explicitly in terms of the operators and connection. This has possible applications to index problems, which are often defined in terms of the Kasparov product. Since Chern character formulas in cyclic homology are most easily computed for unbounded representatives (this is Connes' quantized calculus), explicit representatives of the Kasparov product are desirable in such problems. \newline

Viewing spectral triples as noncommutative metric spaces, the notion of morphism introduced here might shed light on the purely commutative problem of what the correct notion of morphism between metric spaces should be.
\section{Spectral triples and noncommutative geometry}

 By the Gelfand-Naimark theorem, $C^{*}$-algebras can be viewed as noncommutative, locally compact Hausdorff topological spaces. This is the starting point for noncommutative geometry. A continuous map $f:X\rightarrow Y$ between compact Hausdorff spaces gives a *-homomorphism $f^{*}:C(Y)\rightarrow C(X)$ between the dual $C^{*}$-algebras. Noncommutative algebras, however, might not admit any nontrivial algebra homomorphisms. For intstance, the matrix algebra $M_{n}(\C)$ is a such a \emph{simple} algebra. A more flexible notion of morphism for noncommutative algebras is that of a suitable class of bimodules $_{A}\mathpzc{E}_{B}$, with composition coming from the tensor product of bimodules.\newline\newline
Topological $K$-theory is the tool that generalizes in the most straightforward way from spaces to $C^{*}$-algebras. From the definition of $K$-theory it follows readily that $K_{*}(A)\cong K_{*}(M_{n}(A))$ for any $C^{*}$-algebra $A$. This is one of the reasons why one wants to regard these algebras as being equivalent. The notion of \emph{Morita equivalence} formalizes this notion of equivalence and is compatible with the notion of bimodule morphism.\newline\newline
 A Riemannian manifold $M$ is a topological space with some finer structure defined on it. This can be encoded by considering some (pseudo) differential operators on the manifold, e.g. a Dirac operator (when $M$ is Spin$^{c}$), or a signature type operator.\newline\newline
In the spin case, the Riemannian metric on $M$ can be recoverd from the Dirac operator $D$   by
\[d(x,y)=\sup\{\|f(x)-f(y)\|: \|[D,f]\|\leq 1\}. \]
The reader can consult \cite{Conbook} for a proof of this. Recently, some stronger reconstruction theorems have been announced \cite{Conre}.\newline\newline
This motivates the notion of spectral triple \cite{Conspec}.

\begin{definition}A \emph{spectral triple} $(A,\mathpzc{H},D)$ consists of a $\Z/2$-graded $C^{*}$-algebra $A$ represented on a likewise graded Hilbert space $\mathpzc{H}$, together with an odd, selfadjoint operator $D$, with compact resolvent, such that
\[\{a\in A: [D,a]\in B(\mathpzc{H})\},\]
is dense in $A$.
\end{definition}
 Commutative examples are plentiful, mainly given by manifolds. Other examples come from groups, group actions, and foliations. Also, there are various extensions of the notion of spectral triple, notably in the type II and type III setting. Again we refer to \cite{CPR} for these topics.
\section{The noncommutative torus}
The subject of these notes is a notion of \emph{morphism} for spectral triples, a generalization of maps between manifolds. Let us first discuss an example to illustrate this. It will be a noncommutative geometry description of the \emph{fibration} of the torus $S^{1}\times S^{1}$ over the circle $S^{1}$. The projection $S^{1}\times S^{1}\rightarrow S^{1}$ on either of the coordinates is a smooth map, and the fiber over each point is again diffeomorphic to $S^{1}$. Of course this is a very simple fibration because it is just a direct product. However, its noncommutative analogue is very instructive in illustrating the general theory that follows.\newline\newline
The noncommutative torus $A_{\theta}$  is the $C^{*}$-algebra crossed product of the action of $\Z$ on the circle $S^{1}$ by a rotation over the angle $2\pi\theta$, denoted $x\mapsto \alpha_{\theta}(x)$. The algebra $C_{c}(S^{1}\times \Z)$ carries a convolution product
\[f*g(x,n)=\sum_{k\in \Z}f(x,k)g(\alpha_{\theta}^{k}(x),n-k),\]
defining a representation on $\mathpzc{H}:=L^{2}(S^{1}\times \Z)$, yielding the $C^{*}$-algebra $A_{\theta}$.\newline\newline
Another way to describe $A_{\theta}$ is as the universal $C^{*}$-algebra generated by two unitaries $u,v$ subject to the relation $uv=e^{2\pi i\theta} vu$. In this picture, elements of $A_{\theta}$ can be described as series
\[\sum_{n,m\in\Z}\lambda_{n,m}u^{n}v^{m},\]
convergent in a certain norm, analogous to Fourier series.\newline\newline 
The algebra $A_{\theta}$ carries two canonical unbounded derivations, defined on $C_{c}^{\infty}(S^{1}\times\Z)$ by 
\[\partial_{1}f(x,n):= n f(x,n),\quad \partial_{2}f(x,n):=\frac{1}{2\pi i}\partial f(x,n).\]
In the $u,v$ picture, these derivations are 
\[\partial_{1} u^{n}v^{m}=mu^{n}v^{m},\quad \partial_{1} u^{n}v^{m}=nu^{n}v^{m}.\]
On $\mathpzc{H}\oplus\mathpzc{H}$ this yields the operator
\[D:=\begin{pmatrix} 0 & \partial_{1}-i\partial_{2}\\\partial_{1}+i\partial_{2} &0\end{pmatrix}.\]
giving the canonical spectral triple on $A_{\theta}$.

\section{``Fibration" over the circle}
We now describe the structure that we think of as implementing the fibration of $A_{\theta}$ over the circle algebra $C(S^{1})$. It consists of an $(A_{\theta},C(S^{1}))$-bimodule, equipped with an unbounded operator and a connection. The precise structures present on these modules will be described later in these notes. The reader is encouraged to keep this example in mind. \newline\newline
Consider the module $\mathpzc{E}=\ell^{2}(\Z)\tildeotimes C(S^{1})\cong L^{2}(S^{1})\tildeotimes C(S^{1})$. Here $\tildeotimes$ denotes a certain completed tensor product. It 
carries an unbounded, $C(S^{1})$ linear operator
\[S:e_{n}\otimes f\mapsto ne_{n}\otimes f.\]
The canonical spectral triple for the circle $(C(S^{1}),L^{2}(S^{1}),\frac{1}{2\pi i}\partial)$ defines a module of 1-``forms"
\[\Omega^{1}_{\partial}:=\{\sum f_{k}[\frac{1}{2\pi i}\partial,g_{k}]: f_{k},g_{k}\in \textnormal{Lip}^{1}(S^{1})\}\subset B(L^{2}(S^{1})),\]
where Lip$^{1}$ denotes the Lipschitz functions on $S^{1}$.
The module $\mathpzc{E}$ carries a densely defined connection
\[\nabla: e_{n}\otimes f\mapsto e_{n}\otimes[\frac{1}{2\pi i}\partial,f].\] 
$\nabla$ is defined on a dense Lip$^{1}(S^{1})$-submodule $E^{1}\subset\mathpzc{E}$, and maps it into $E^{1}\otimes_{\textnormal{Lip}^{1}}\Omega^{1}_{\partial}$.
It satisfies $[\nabla,S]=0$. The tensor product $\mathpzc{E}\tildeotimes_{C(S^{1})} L^{2}(S^{1})$ is isomorphic to $\mathpzc{H}$. Under this identification, the derivation $\partial_{2}$ equals 
\[e\otimes h\mapsto e\otimes \frac{1}{2\pi i}\partial h + \nabla(e)h.\]
This expression is well defined because $\nabla$ satisfies a Leibniz rule
\[\nabla(ef)=\nabla(e)f+ e\otimes [\frac{1}{2\pi i}\partial,f] .\]
We denote it by $1\otimes_{\nabla}\frac{1}{2\pi i}\partial$. We thus see that the canonical spectral triple on $A_{\theta}$ can be factorized as a graded tensor product
\[(A_{\theta},\mathpzc{H}\oplus\mathpzc{H},D)=(E,S,\nabla)\otimes(C(S^{1}),L^{2}(S^{1}),\frac{1}{2\pi i}\partial).\]
The tensor product on the right is to be interpreted as
\[\mathpzc{E}\otimes_{C(S^{1})} L^{2}(S^{1})\oplus\mathpzc{E}\otimes_{C(S^{1})} L^{2}(S^{1}),\]
with operator
\[\begin{pmatrix} 0& S\otimes 1-i1\otimes_{\nabla}\frac{1}{2\pi i} \partial\\
S\otimes 1 + i1\otimes_{\nabla}\frac{1}{2\pi i} \partial &0\end{pmatrix}.\]

 Thus, by choosing the right gradings, the triple $(\mathpzc{E},S,\nabla)$ can be viewed as a  fibration of the noncommutative torus over the circle.

\section{$C^{*}$-modules and regular operators}
We now proceed by describing the modules, operators and connections involved in a more rigorous manner. Let $(A,B)$ be a pair of separable, $\Z/2$-graded $C^{*}$-algebras. The reader who feels uneasy thinking about graded $C^{*}$-algebras, can think of trivially graded (i.e. ungraded) $C^{*}$-algebras. The reason for developing the theory for graded algebras is that one can treat the even and odd cases of $K$-theory at the same time. The standard reference for the theory of $C^{*}$-modules is \cite{Lan}.
\begin{definition}A \emph{$C^{*}$-module over $B$} is a right $B$-module $\mathpzc{E}$ equipped with a positive definite $B$-valued inner product.
\end{definition}
A positive definite $B$-valued inner product is a pairing $\mathpzc{E}\times\mathpzc{E}\rightarrow B$, satisfying
\begin{itemize}\item $\langle e_{1},e_{2}\rangle=\langle e_{2},e_{1}\rangle^{*},$
\item $\langle e_{1},e_{2}b\rangle=\langle e_{1},e_{2}\rangle b,$ 
\item $\langle e,e\rangle\geq 0$ and $\langle e,e\rangle=0\Leftrightarrow e=0,$\item $\mathpzc{E}$ is complete in the norm $\|e\|^{2}:=\|\langle e,e\rangle\|.$\end{itemize}
We use the notation $\mathpzc{E}\leftrightharpoons B$ to indicate this structure. \newline\newline
The natural endomorphisms to consider in a $C^{*}$-module are the following:
\[\Endst_{B}(\mathpzc{E}):=\{T:\mathpzc{E}\rightarrow \mathpzc{E}:\exists T^{*}:\mathpzc{E}\rightarrow\mathpzc{E}, \langle Te,f\rangle=\langle e, T^{*}f\rangle\}.\] 
Operators in $\Endst_{B}(\mathpzc{E})$ are automatically $B$-linear and bounded, and they form a $C^{*}$-algebra in the operator norm and the involution $T\mapsto T^{*}$.\newline\newline
There is a natural $C^{*}$-subalgebra, analogous to the compact operators on a Hilbert space.  The algebra of \emph{compact endomorphisms} $\K_{B}(\mathpzc{E})\subset \Endst_{B}(\mathpzc{E})$ is the $C^{*}$-subalgebra generated by the operators $e\otimes f (g):=e\langle f,g\rangle$.\newline\newline
An $(A,B)$-\emph{bimodule} is a $C^{*}$-module $\mathpzc{E}\lrh B$, together with a graded *-homomorphism $A\rightarrow\Endst_{B}(\mathpzc{E})$.\newline\newline
A \emph{regular operator in} $\mathpzc{E}$ is a densely defined  closed operator, $D:\Dom D\rightarrow \mathpzc{E}$, such that $D^{*}$ is densely defined in $\mathpzc{E}$ and $1+D^{*}D$ has dense range. This condition is automatic in the Hilbert space setting, but needs to be imposed in $C^{*}$-modules, to avoid pathologies. 
The operator $D$ is \emph{selfadjoint} if it is symmetric on its domain, and $\Dom D^{*}=\Dom D$. An excellent reference for the theory of regular operators in $C^{*}$-modules is \cite{Lan}.

\section{Unbounded $KK$-theory}
$KK$-theory associates to a pair $(A,B)$ of separable, $\Z/2$-graded  $C^{*}$-algebras a $\Z/2$-graded abelian group $KK_{*}(A,B)$. Kasparov \cite{Kas} originally constructed and described these groups using bounded \emph{Fredholm operators} in $C^{*}$-modules.\newline\newline
A defining element of the group $KK_{0}(A,B)$ is a pair $(\mathpzc{E},F)$ consisting of an $(A,B)$-bimodule $\mathpzc{E}$, together with an operator $F\in\Endst_{B}(\mathpzc{E})$ satisfying
\[a(F^{2}-1),\quad a(F-F^{*}),\quad [F,a] \in \K_{B}(\mathpzc{E}).\]
Subsequently one considers unitary equivalence classes of such pairs, and quotients by the relation of homotopy to obtain the abelian group $KK_{0}(A,B)$. The groups $KK_{i}(A,B)$ are defined as being $KK_{0}(A, B\tildeotimes \C_{i})$, where $\C_{i}$ is the $i$-th complex Clifford algebra. This is a graded $C^{*}$-algebra, and it is at this point that working with graded algebras comes in handy.
\newline\newline
Kasparov's main achievement was the construction of an associative, distributive product
\[KK_{i}(A,B)\otimes_{\Z} KK_{j}(B,C)\rightarrow KK_{i+j}(A,C),\]
now known as the \emph{Kasparov product}.
The Kasparov product has remarkable properties. It allows one to view the $KK$-groups as the morphisms in a category $\mathfrak{KK}$ whose objects are $C^{*}$-algebras. Moreover, Cuntz \cite{Cuntz} and Higson \cite{Hig} showed that that $KK$-theory has a \emph{universal property}, in the sense that any functor from $C^{*}$-algebras to abelian groups which is Morita invariant and split exact, factors through this category $\mathfrak{KK}$. Such functors are automatically homotopy invariant. In this sense $KK$-theory is the universal cohomology theory for $C^{*}$-algebras.\newline

In the above Fredholm picture, the Kasparov product is very difficult to define, and we will refrain from doing so here. We will describe the product in a different picture, given below.

\begin{definition}[\cite{BJ}]The cycles for $KK_{0}(A,B)$ may also be described by pairs $(\mathpzc{E},D)$, where
\begin{itemize}

\item $\mathpzc{E}$ is an $(A,B)$-bimodule.
\item$D:\Dom D\rightarrow \mathpzc{E}$ is an odd selfadjoint regular operator.
\item$\forall a\in A: a(1+D^{2})^{-1}\in\K_{B}(\mathpzc{E})$.
\item The subalgebra
\[\mathcal{A}_{1}:=\{a\in A:[D,a]\in\Endst_{B}(\mathpzc{E})\},\]
is dense in $A$.
\end{itemize}
Such pairs $(\mathpzc{E},D)$ are referred to as \emph{unbounded $KK$-cycles}.
\end{definition}
The relation between the bounded and the unbounded picture is given by a simple procedure.
The following results are due to Baaj-Julg\cite{BJ}. 
\begin{theorem}[\cite{BJ}] Let $F:=D(1+D^{2})^{-\frac{1}{2}}\in\Endst_{B}(\mathpzc{E})$, the \emph{bounded transform} of $D$.
\begin{itemize}
\item $(\mathpzc{E},F)$ is a \emph{Kasparov module}, i.e. $F^{*}=F$ and 
\[ \forall a\in A, a(F^{2}-1), [F,a]\in\K_{B}(\mathpzc{E}).\]
\item Two unbounded modules are \emph{equivalent} if their bounded transforms are homotopic. Any Kasparov module is homotopic to the bounded transform of an unbounded one.
\end{itemize}
\end{theorem}
Their motivation for introducing the unbounded picture was that it simplifies another product structure in Kasparov's theory, the \emph{external product}
\[KK_{i}(A,B)\otimes KK_{j}(A',B')\rightarrow KK_{i+j}(A\minotimes A',B\minotimes B'),\] 
where $A,A',B,B'$ are distinct $C^{*}$-algebras. Baaj and Julg proved the following
\begin{theorem}[\cite{BJ}]

On unbounded cycles, the \emph{external Kasparov product} 
is given by
\[(\mathpzc{E},S)\times(\mathpzc{F},T):=(\mathpzc{E}\minotimes\mathpzc{F}, S\otimes 1 + 1\otimes T),\]
where
\[1\otimes T(e\otimes f):=(-1)^{\partial e}e\otimes Tf.\]
\end{theorem} 
In the case $B=B'=\C$, this product corresponds to the direct product of manifolds. The case $A=A'=\C$ gives the external product in topological $K$-theory.
\section{Algebraic intermezzo}
When trying to define the \emph{internal Kasparov product}
\[KK_{i}(A,B)\otimes KK_{j}(B,C)\rightarrow KK_{i+j}(A,C),\] on unbounded cycles, we run into the following problem. In the Fredholm picture, Kasparov proved that on the module $\mathpzc{E}\tildeotimes_{B}\mathpzc{F}$ one can always find an operator, unique up to homotopy, that defines the class of the Kasparov product. In the unbounded picture, as in the case of the external product, the natural guess for the operator is something like $S\otimes 1+ 1\otimes T$. However, the expression $1\otimes T$ does not make sense, since $T$ does not commute with the elements of $B$, and we take a balanced tensor product. It turns out that there is a notion of \emph{connection} which corrects for this problem. The algebraic theory of forms and connections is described in detail in \cite{CQ}. \newline\newline
For clarity, we first consider the following structure of a category on algebraic $(A,B)$-bimodules with odd operator $(E,D)$.
\begin{definition} Let $B$ be an algebra. The module of $1$-\emph{forms} of $B$ is the kernel of the graded multiplication map
\[\begin{split}\Omega^{1}(B):=\ker(B\otimes B&\xrightarrow{m} B)\\
b_{1}\otimes b_{2} &\mapsto b_{1}\gamma(b_{2}),\end{split}\]
where $\gamma\in\Aut B$ is the grading automorphism.
The \emph{universal derivation} $d:B\rightarrow \Omega^{1}(B)$ is given by
\[b\mapsto 1\otimes b- \gamma(b)\otimes 1.\]
\end{definition}
Any derivation $\delta:B\rightarrow M$ into a $B$-bimodule $M$ factors through the bimodule $\Omega^{1}(B)$ in the following sense.
\begin{proposition}[\cite{CQ}] The bimodule $\Omega^{1}(B)$ is universal for derivations $\delta:B\rightarrow M$, where $M$ is a $B$-bimodule. That is, for any such $\delta$ there is a unique map $j_{\delta}:\Omega^{1}(B)\rightarrow M$ such that $\delta=j_{\delta}\circ d$.
\end{proposition}
The map $j_{\delta}$ is defined by setting $j_{\delta}(da)=\delta(a)$. This determines $j_{\delta}$ as a bimodule map, because the elements $da$ generate $\Omega^{1}(B)$ as a bimodule.

\begin{definition}
A \emph{connection} on a right $B$-module $E$ is a map \[\nabla:E\rightarrow E\otimes_{B}\Omega^{1}(B),\] satisfying
\[\nabla(eb)=\nabla(e)b+e\otimes db.\]
\end{definition}
If a connection $\nabla$ on $E$ is given, $F$ is a $(B,C)$-bimodule and $T\in \End_{B}(F)$, then the operator
\[1\otimes_{\nabla}T(e\otimes f):=(-1)^{\partial e\partial T}(e\otimes Tf+\nabla_{T}(e)f),\]
is well defined on $E\otimes_{B} F$. Here $\partial e,\partial T\in\{0,1\}$ denote the degree of the homogeneous elements $e$ and $T$ respectively. The connection $\nabla_{T}:E\rightarrow E\otimes_{B}\End_{C}(F)$ is the composition $j_{\delta}\circ\nabla$ with $\delta$ the derivation $b\mapsto [T,b]$. When a connection $\nabla '$ is given on $F$, we can apply the same trick and define a connection
\[1\otimes_{\nabla}\nabla ':E\otimes_{B}F\rightarrow E\otimes_{B} F\otimes_{C}\Omega^{1}(C),\]
now by using the derivation $b\mapsto [\nabla,b].$ An \emph{isomorphism} of triples $(E,S,\nabla)$ and $(E',S',\nabla ')$ is a bimodule isomorphism $g:E\rightarrow F$ with the additional properties that
\begin{itemize}
\item $g^{-1}S'g=S$;
\item $g^{-1}\nabla' g=\nabla$.
\end{itemize}
Of course, isomorphism of triples is an equivalence relation.

\begin{proposition}[\cite{Mes}]\label{alg} Let $A,B,C$ be algebras, $E,F$ $(A,B)-$ and $(B,C)-$ bimodules respectively. The composition law
\[(E,S,\nabla)\circ (F,T,\nabla '):=(E\otimes_{B} F, S\otimes 1+ 1\otimes_{\nabla}T, 1\otimes_{\nabla}\nabla '),\]
is associative up to isomorphism. Isomorphism classes of triples $(E,S,\nabla)$ are the morphisms in a category whose objects are pairs $(E,D)$, where $E$ is an $(A,B)$-bimodule and $D\in \End_{B}(E)$ an endomorphism.
\end{proposition}
\begin{remark} A morphism from $(G,D)$ to $(F,T)$ is a triple $(E,S,\nabla)$ such that $(E\otimes_{B} F, S\otimes 1+1\otimes_{\nabla}T)$ is isomorphic to $(G,D)$.
\end{remark}
In this setting a spectral triple $(A,\mathpzc{H},D)$ is more conveniently denoted by just $(\mathpzc{H},D)$. In particular, these are $(A,\C)$ bimodules. Unfortunately, the algebraic setting discussed above is not appropriate for dealing with spectral triples. It needs to be enriched to accommodate for the analytic phenomena governing them.

In order to construct a category of spectral triples (or unbounded bimodules) in which the morphisms are unbounded bimodues $(\mathpzc{E},D)$, with some notion of connection, several problems need to be addressed:
\begin{itemize}
\item Unbounded regular operators are not endomorphisms (i.e. not everywhere defined).
\item The graded commutators $[D,a]$ are endomorphisms only for $a$ in a dense subalgebra of $A$.
\item An analytic version of $\Omega^{1}(B)$ and the notion of connection for dense subalgebras are needed.
\item The product operator $S\otimes 1+ 1\otimes_{\nabla}T$ should be selfadjoint, regular and have compact resolvent.
\end{itemize}
All these issues can be resolved by introducing an appropriate notion of smoothness for unbounded $KK$-cycles.

\section{Operator algebras and modules}
To overcome the aforementioned problems, we need to broaden our scope from $C^{*}$-algebras to operator spaces. The algebraic structures of algebras and modules will need operator space analogues as well. Operator space theory was developed by Effros and Ruan \cite{EfR},\cite{Ruan}, and many others.
\begin{definition}
An \emph{operator space} is a closed subspace of some $C^{*}$-algebra. 
\end{definition}
The main feature of an operator space $X$ is that it comes with canonical matrix norms, i.e. $M_{n}(X)$ carries a canonical norm.
A map $\phi:X\rightarrow Y$ between operator spaces is \emph{completely bounded} if $\|\phi\|_{cb}:=\sup_{n}\|\phi_{n}\|<\infty$, where $\phi_{n}:M_{n}(X)\rightarrow M_{n}(Y)$ is the map induced by $\phi$. It is \emph{completely contractive} if $\|\phi\|_{cb}\leq 1$. The completely bounded maps form the natural class of maps between operator spaces.
\begin{example} A *-homomorphism $\phi:A\rightarrow B$ between $C^{*}$-algebras is automatically completely bounded, as is an adjointable operator $T\in\Endst_{B}(\mathpzc{E},\mathpzc{F})$ between $C^{*}$-modules.
\end{example}

The natural tensor product for operator spaces $X$ and $Y$ is the \emph{Haagerup tensor product}, denoted by $X\tildeotimes Y$. Its norm is given by
\[\|z\|:=\inf\{\|\sum x_{i}x_{i}^{*}\|^{\frac{1}{2}}\|\sum y_{i}^{*}y_{i}\|^{\frac{1}{2}}:z=\sum x_{i}\otimes y_{i}\}.\]
Note that although $x^{*}$ need not be an element of $X$, it does make sense in the containing $C^{*}$-algebra of $X$. The space $X\tildeotimes Y$ is again an operator space.
 An \emph{operator algebra} is an operator space $\mathcal{A}$ whose multiplication $\mathcal{A}\tildeotimes \mathcal{A}\rightarrow\mathcal{A}$ is completely contractive.
An \emph{involutive operator algebra} is an operator algebra with an involution $a\mapsto a^{*}$ which is completely bounded.
 An \emph{operator module} $M$ over an operator algebra $\mathcal{B}$ is an operator space $M$, which is also a (say) right $\mathcal{B}$-module, such that the module action $M\tildeotimes\mathcal{B}\rightarrow M$ is completely bounded.
 The \emph{Haagerup module tensor product} $M\tildeotimes_{\mathcal{B}}N$ of right and left $\mathcal{B}$ operator modules $M$ and $N$, respectively, is the quotient of $M\tildeotimes N$ by the closed subspace generated by $mb\otimes n-m\otimes bn$. The reader can consult \cite{Blechbook} and \cite{BMP} for many aspects of the theory of operator modules. Also, see \cite{Helem} for a survey on operator space tensor products. 
\begin{example} A $C^{*}$-module $\mathpzc{E}\lrh B$ is canonically an operator space by viewing it as the upper right corner of its \emph{linking algebra} $\K_{B}(B\oplus\mathpzc{E})$. As such it is an operator module and the Haagerup tensor product of $C^{*}$-modules is completely isometrically isomorphic to the $C^{*}$-module tensor product.
\end{example}
D.P. Blecher \cite{Blech} observed that, when $\mathpzc{E}$ is countably generated, by choosing an approximate unit
\[u_{n}=\sum_{1\leq |i|\leq n} x_{i}\otimes x_{i}\in\K_{B}(\mathpzc{E}), \]
(which is possible by Kasparov's stabilization theorem \cite{Kas})
$\mathpzc{E}$ can be written as an inductive limit of the canonical modules $B^{2n}.$
This is done by considering the maps 
\[\phi_{n}:(b_{i})\mapsto\sum_{1\leq |i|\leq n} x_{i}b_{i},\quad\psi_{n}:e\mapsto (\langle x_{i},e\rangle),\]
which are completely contractive, and $\phi_{n}\circ\psi_{n}\rightarrow 1$ strongly. The inner product in $\mathpzc{E}$ can be recovered from these maps as
\[\langle e,f\rangle:=\lim_{\alpha}\langle \psi_{n}(e),\psi_{n}(f)\rangle_{n}.\]

\section{Stably rigged modules}
Blecher used his observation to develop a theory of modules over operator algebras, that are in many ways similar to $C^{*}$-modules. In case the algebra is actually a $C^{*}$-algebra, this class of modules coincides with that of $C^{*}$-modules. See \cite{Blech2} for details.
\begin{definition}[Blecher] Let $\mathcal{B}$ be an operator algebra with contractive countable approximate identity. A \emph{rigged module} over $\mathcal{B}$ is a right operator module $E$ over $\mathcal{B}$ together with completely contractive module maps $\psi_{n}:E\rightarrow \mathcal{B}^{2n}$ and $\phi_{n}:\mathcal{B}^{2n}\rightarrow E$, such that $\phi_{n}\circ\psi_{n}$ converges strongly to $1$, and $\phi_{n}$ is $\mathcal{B}$-essential. When $\phi_{n},\psi_{n}$ and the approximate identity are merely completely bounded, $E$ is an \emph{stably rigged module}.
\end{definition}
The difference between rigged and stably rigged modules might seem only formal at first sight. However, the contractivity assumption is a fairly strong one. 

\begin{theorem}[Blecher] A rigged module over a $C^{*}$-algebra is a $C^{*}$-module, and the Haagerup tensor product of (stably) rigged modules is again a (stably) rigged module.\end{theorem}

A rigged module over a $C^{*}$-algebra is completely isometrically isomorphic to a $C^{*}$-module. In the cb-setting such a theorem has not been established, and can definitely not be proven in a similar way. An important corollary of the above theorem is that for a rigged module $E$ over an operator algebra $\mathcal{B},$ and a completely contractive homomorphism $\mathcal{B}\rightarrow\Endst_{C}(\mathpzc{F})$, with $\mathpzc{F}\lrh C$ a $C^{*}$-module, the Haagerup tensor product $E\tildeotimes_{\mathcal{B}}\mathpzc{F}$ is a genuine $C^{*}$-module. This fact will be exploited when dealing with graphs of unbounded operators.



\section{Sobolev modules}
This section describes the construction of Sobolev modules and algebras as developed in \cite{Mes}. They are the analogues of the usual Sobolev spaces that appear in Riemannian geometry, but we describe them in a more algebraic manner. For this reason we obtain only Sobolev spaces indexed by the natural numbers, as opposed to the positive real numbers. This is to avoid the use of functional calculus.\newline\newline

The \emph{graph} of a regular operator $D$ in $\mathpzc{E}$ is the closed submodule
\[\mathfrak{G}(D):=\{(e,De):e\in\Dom D\}\subset \mathpzc{E}\oplus\mathpzc{E}.\]
If $D$ is selfadjoint, we define
\[\Dom D_{2}:=\{(e,De)\in\mathfrak{G}(D):e\in\Dom D^{2}\},\]
and
\[D_{2}:(e,De)\mapsto (De,D^{2}e).\]
\begin{lemma} Let $D$ be a selfadjoint regular operator in $\mathpzc{E}$. The operator $D_{2}$ is selfadjoint and regular in $\mathfrak{G}(D)$.
\end{lemma}

Iterating this construction gives the \emph{Sobolev chain} of $D$:
\[\cdots\rightarrow\mathfrak{G}(D_{n+1})\rightarrow\mathfrak{G}(D_{n})\rightarrow\cdots\rightarrow\mathfrak{G}(D_{2})\rightarrow\mathfrak{G}(D)\rightarrow\mathpzc{E}.\]
In $C^{*}$-modules, not every closed submodule is the range of a projection in $\Endst_{B}(\mathpzc{E})$. Modules with this property are called \emph{complemented submodules}. The following theorem states that the graph of a regular operator in a $C^{*}$-module is a complemented submodule. The regularity condition on unbounded operators is imposed mainly for this reason.
\begin{theorem}[\cite{Baaj},\cite{Lan},\cite{Wor}] Let $D$ be a selfadjoint regular operator in $\mathpzc{E}$. Then
\[p_{D}:=\begin{pmatrix}(1+D^{2})^{-1} & D(1+D^{2})^{-1}\\D(1+D^{2})^{-1} & D^{2}(1+D^{2})^{-1}\end{pmatrix},\]
is a projection in $\Endst_{B}(\mathpzc{E}\oplus\mathpzc{E})$, and $p(\mathpzc{E}\oplus\mathpzc{E})=\mathfrak{G}(D)$. Moreover

\[\mathfrak{G}(D)\oplus v\mathfrak{G}(D)\cong\mathpzc{E}\oplus\mathpzc{E},\]
is an orthogonal direct sum, where $v$ is the unitary $v:(x,y)\mapsto(-y,x)$.
\end{theorem}
This result is attributed to several people, but Woronowicz explicitly mentions the projection $p_{D}$, which is why we refer to it as the Woronowicz projection. The Sobolev modules and Woronowicz projections can be used to construct a chain of subalgebras
\[\cdots \subset\mathcal{A}_{k+1}\subset\mathcal{A}_{k}\subset\cdots\subset \mathcal{A}_{1}\subset A,\]
for any spectral triple or $KK$-cycle, in the following way.
For a $KK$-cycle $(\mathpzc{E},D)$, we have a representation
\[\begin{split}\pi_{1}:\mathcal{A}_{1} &\rightarrow\Endst_{B}(\mathpzc{E}\oplus\mathpzc{E})\cong M_{2}(\Endst_{B}(\mathpzc{E}))\\
a &\mapsto\begin{pmatrix}a & 0\\ [D,a] & (-1)^{\partial a} a\end{pmatrix}.\end{split}\]
This gives a representation
\[\begin{split}\theta_{1}:\mathcal{A}_{1}&\rightarrow\Endst_{B}(\mathfrak{G}(D))\oplus\Endst_{B}(v\mathfrak{G}(D))\\
a&\mapsto p_{D}\pi_{1}(a)p_{D}+p^{\perp}_{D}\pi_{1}(a)p^{\perp}_{D}.\end{split}\]
The restriction $\chi_{1}$ of $\theta_{1}$ to $\mathfrak{G}(D)$ acts as
\[\chi_{1}(a):\begin{pmatrix}e\\De\end{pmatrix}\mapsto\begin{pmatrix}ae\\Dae\end{pmatrix}.\]

This allows us to inductively define
\[\mathcal{A}_{n+1}:=\{a\in\mathcal{A}_{n}:[D,\theta_{n}(a)]\in\Endst_{B}(\mathfrak{G}(D_{n}))\},\]
\begin{equation}\label{pi}\begin{split}\pi_{n+1}:\mathcal{A}_{n}&\rightarrow\Endst_{B}(\mathfrak{G}(D_{n})\oplus\mathfrak{G}(D_{n}))\\
a&\mapsto\begin{pmatrix}\theta_{n}(a) & 0\\
[D,\theta_{n}(a)] & (-1)^{\partial a}\theta_{n}(a)\end{pmatrix},\end{split}\end{equation}

\begin{equation}\label{the}\theta_{n+1}(a):=p_{n+1}p_{n}\pi_{n+1}(a)p_{n}p_{n+1}+p_{n+1}^{\perp}p_{n}^{\perp}\pi_{n+1}(a)p_{n}^{\perp}p_{n+1}^{\perp}.\end{equation}

\begin{definition}The algebra $\mathcal{A}_{n}$ is the $n$-th \emph{Sobolev subalgebra} of $A$. It allows for a completely contractive representation $\chi_{n}:\mathcal{A}_{n}\rightarrow\mathfrak{G}(D_{n})$, which is \emph{not} a *-homomorphism. When $A=\Endst_{B}(\mathpzc{E})$, we write $\Sob_{n}(D)$ for $\mathcal{A}_{n}$.
\end{definition}
The algebras $\mathcal{A}_{n}$ can also be characterized by a relative boundedness condition.
\begin{proposition}[\cite{Mes}] \label{rel}We have $a\in\mathcal{A}_{n}$ if and only if \[(\textnormal{ad} (D))^{n}(a)(D\pm i)^{-n+1},(\textnormal{ad} (D))^{n}(a^{*})(D\pm i)^{-n+1}\in\Endst_{B}(\mathpzc{E}).\]
\end{proposition}
The representations $\bigoplus_{j=0}^{n}\pi_{j}$ realize $\mathcal{A}_{n}$ a closed subspace of a $C^{*}$-algebra, i.e. as an \emph{operator space}. Taking these representations as defining the topology on $\mathcal{A}_{n}$, the inclusions $\mathcal{A}_{n+1}\rightarrow\mathcal{A}_{n}$ become completely contractive *-homomorphisms.
\begin{proposition}[\cite{Mes}] The involution on $\mathcal{A}_{n}$ is a complete anti-isometry.
\end{proposition}
Thus, the Sobolev subalgebras are involutive operator algebras in their natural operator space topology.
\section{Smoothness}
Although the Sobolev subalgebras of a given $KK$-cycle always exist and contain the identity, in general we know very little about them. One of the conditions in the definition of $KK$-cycle is that the algebra $\mathcal{A}_{1}$ is dense in the $C^{*}$-algebra $A$. This can be interpreted as a smoothness condition.
\begin{definition} A $KK$-cycle $(E,D)$ is said to be \emph{(left) $C^{k}$} if $\mathcal{A}_{k}$ is dense in $A$. It is said to be \emph{(left) smooth} is it is (left) $C^{k}$ for all $k$ and $\mathscr{A}=\bigcap_{k}\mathcal{A}_{k}$ is dense in $A$.
\end{definition}
This definition of smoothness is weaker then the one employed in \cite{Conspec}. In particular, spectral triples coming from manifolds are smooth in our sense.
Indeed, for a $C^{k}$-cycle the Sobolev algebras have good properties.
\begin{theorem}[\cite{Mes}]  If $(E,D)$ is $C^{k}$ then the algebras $\mathcal{A}_{i}$ for $i\leq k$ are stable under holomorphic functional calculus in $A$.
\end{theorem}
\begin{definition} A \emph{smooth} $C^{*}$-\emph{algebra} is an inverse system of involutive operator algebras
\[\cdots\rightarrow \mathcal{A}_{n+1}\rightarrow\mathcal{A}_{n}\rightarrow\cdots\rightarrow A,\]
coming from a spectral triple.
\end{definition}
Smooth $C^{*}$-algebras should be thought of as the analogues of smooth manifolds. By holomorphic stability, any finitely generated projective module over a unital smooth $C^{*}$-algebra can be smoothened. In the case of countably generated modules, smoothness is not direct anymore, and needs to be imposed on the module.
\begin{definition} Let $B$ be a $C^{k}$-algebra, and $\mathpzc{E}\lrh B$ a $C^{*}$-module. $\mathpzc{E}$ is said to be $C^{k}$ if there is an approximate unit
\[u_{n}:=\sum_{1\leq |i|\leq n}x_{i}\otimes x_{i}\in \K_{B}(\mathpzc{E}),\]
such that the norm of the infinite matrix
\[\|(\langle x_{i},x_{j}\rangle)\|_{k}\leq C.\]

\end{definition}
\begin{remark} The $k$-norm is the norm induced by the representation $\bigoplus_{j=0}^{k}\pi_{j}$. Since this is an operator norm, it gives norms for all matrix algebras.
\end{remark}
The approximate unit, the existence of which is demanded, can be used to construct a chain of submodules
\[\cdots\subset E^{k+1}\subset E^{k}\subset\cdots\subset E^{1}\subset\mathpzc{E},\]
which correspond to higher order Lipschitz sections of a vector bundle. In the finitely generated unital case the approximate unit is an actual unit. It is no more than a choice of projection in the subalgebra, which, as mentioned above, is always possible.

\begin{proposition}\label{useful} Let $B$ be a smooth $C^{*}$-algebra and $\mathpzc{E}$ a  $C^{k}$-$B$-module. Then
\[E^{k}:=\{e\in\mathpzc{E}:\langle x_{i},e\rangle\in\mathcal{B}_{k},\quad\sup_{n}\|\sum_{1\leq |i|\leq n}e_{i}\langle x^{\alpha}_{i},e\rangle\|_{k}<\infty\},\]
is an stably rigged $\mathcal{B}_{k}$-module. Moreover, the inclusions $E^{k+1}\rightarrow E^{k}$ are completely contractive with dense range, and $E^{k+1}\tildeotimes_{\mathcal{B}_{k+1}}\mathcal{B}_{k}\cong E^{k}$.
\end{proposition}
The $E^{k}$ are stably rigged modules, but they are constructed in a very specific way. This allows us to say a lot more about them than for general stably rigged modules.
\begin{theorem} Let $\mathpzc{E}$ be a countably generated $C^{k}$-module over a $C^{k}$-algebra $B$. For all $i\leq k$, there are cb-isomorphisms $E^{i}\oplus\mathpzc{H}_{\mathcal{B}_{i}}\cong \mathpzc{H}_{\mathcal{B}_{i}}$, compatible with the $C^{k}$-structure. Consequently, a countably generated $C^{*}$-module is a $C^{k}$ module if and only if it is completely isomorphic to a direct summand in a rigged module.
\end{theorem}
For stably rigged modules, operator algebras $\Endst_{\mathcal{B}}(E)$ and $\K_{\mathcal{B}}(E)$ are defined \cite{Blech2}. For $C^{k}$-modules over a $C^{k}$-algebra, the definitions are the same as in the $C^{*}$-case.
\begin{theorem}[\cite{Mes}] The submodules $E^{i}\subset\mathpzc{E}$, $i\leq k$, inherit a $\mathcal{B}_{i}$-valued inner product by restriction of the inner product on $\mathpzc{E}$. We have
\[\Endst_{\mathcal{B}_{i}}(E^{i})=\{T:E^{i}\rightarrow E^{i}\quad :\quad\exists T^{*}:E^{i}\rightarrow E^{i},\quad\langle Te,f\rangle=\langle e, T^{*}f\rangle\},\]
and $\K_{\mathcal{B_{i}}}({E^{i}})$ is the $i$-operator norm closure of the finite rank operators in $\Endst_{\mathcal{B}_{i}}(E^{i})$. Moreover there is a cb-isomorphism 
\[\K_{\mathcal{B}_{i}}(E^{i})\cong E^{i}\tildeotimes_{\mathcal{B}_{i}} E^{i *},\]
and
\[\K_{\mathcal{B}_{i}}(E^{i})=\K_{B}(\mathpzc{E})\cap\Endst_{\mathcal{B}_{i}}(E^{i}).\]
\end{theorem}
 That is, they are those operators $T:E^{i}\rightarrow E^{i}$ that admit an adjoint with respect to the above inner product. Such operators are automatically completely bounded. The involution $T\mapsto T^{*}$ is in general not a complete isometry, but we have $\frac{1}{C}\|T\|_{i}\leq\|T^{*}\|_{i}\leq C\|T\|_{i}$ for some $C\geq 1$. In particular, unitaries are not necessarily isometries.
\begin{remark}Note that the topology on the $E^{i}$ is not defined by the inner product, but by the approximate unit. 
\end{remark} 
In case we have two smooth $C^{*}$-algebras $A$ and $B$, we can now state what it means for a smooth module to be smooth as a bimodule.

\begin{definition} Let $A,B$ be $C^{k}$-algebras. A $C^{k}$-module $\mathpzc{E}\lrh B$ is a $C^{k}$ bimodule if the $A$ representation restricts to representations $\mathcal{A}_{i}\rightarrow\Endst_{\mathcal{B}_{i}}(E^{i})$, for $i\leq k$.
\end{definition}

\section{Transverse operators}
The theory of regular operators can be developed for $C^{k}$-modules. Definitions and most of the essential results still hold true, but their proofs are quite different from the $C^{*}$ setting. Thus, a selfadjoint operator $D:\Dom D\rightarrow E^{i}$ is said to be regular if it is closed, its domain $\Dom D\subset E^{i}$ is dense in $E^{i},$ and equals the domain of its adjoint and the range of the operators $D\pm i$ is all of $E^{i}$. A selfadjoint regular operator in $E^{k}$ extends to a regular operator in $E^{i}$, $i\leq k$, as $D\otimes 1$, by proposition \ref{useful}. The main result on selfadjoint regular operators in $C^{k}$-modules is the existence of the Woronowicz projection.
\begin{theorem} Let $\mathpzc{E}$ be a $C^{k}$-module over a $C^{k}$-algebra $B$, and $D$ a selfadjoint regular operator in $E^{k}$. Then
\[p_{D}:=\begin{pmatrix}(1+D^{2})^{-1} & D(1+D^{2})^{-1}\\D(1+D^{2})^{-1} & D^{2}(1+D^{2})^{-1}\end{pmatrix},\]
is a projection in $\Endst_{\mathcal{B}_{k}}(E^{k}\oplus E^{k})$, and $p(E^{k}\oplus E^{k})=\mathfrak{G}(D)$. Moreover

\[\mathfrak{G}(D)\oplus v\mathfrak{G}(D)\cong(E^{k}\oplus E^{k}),\]
is an orthogonal direct sum, where $v$ is the unitary $v:(x,y)\mapsto(-y,x)$.
\end{theorem}
This implies that we get Sobolev subalgebras $\Sob_{i}^{k}(D)\subset\Endst_{\mathcal{B}_{k}}(E^{k})$ for all $i$. We can use the same formulae \eqref{pi},\eqref{the} to define the representations $\pi^{k}_{i},\theta^{k}_{i}$of these Sobolev algebras. We get the same relative boundedness conditions as in proposition \ref{rel}, but now for the $i$-norms.

\begin{definition} A $KK$-cycle $(\mathpzc{E},D)$ over $C^{k}$-algebras $(A,B)$ is said to be $C^{k}$ if $\mathpzc{E}$ is a $C^{k}$-bimodule, and $D$ restricts to a regular operator in $E^{k-1}$.
\end{definition}
$D$ is said to be \emph{transverse} $C^{k}$ if $\mathcal{A}_{n}\rightarrow\Sob^{i}_{j}(D)$ completely boundedly, for all $i+j=n\leq k$ (transversality).

\section{Smooth connections}
The Haagerup tensor product linearizes the multiplication in an operator algebra continuously. Since the definition of connections and 1-forms in section 6 essentially only uses the multiplication in an algebra, these definitions carry over to operator algebras. We define \[\Omega^{1}(\mathcal{B}):=\ker(\mathcal{B}\tildeotimes \mathcal{B}\rightarrow \mathcal{B}).\] Connections are defined as in the algebraic setting:
A $C^{k}$-connection in a $C^{k}$-module $\mathpzc{E}$ over a $C^{k}$-algebra $B$ is a connection
\[\nabla:E^{k}\rightarrow E^{k}\tildeotimes_{\mathcal{B}_{k}}\Omega^{1}(\mathcal{B}_{k}),\]
which is completely bounded for the present operator space topologies. Since our modules carry inner products, we now require the extra condition of being a $*$-\emph{connection}. This means there is another connection
\[\nabla^{*}:E^{k}\rightarrow E^{k}\tildeotimes_{\mathcal{B}_{k}}\Omega^{1}(\mathcal{B}_{k}),\]
such that 
\[\langle e_{1},\nabla(e_{2})\rangle - \langle \nabla^{*}(e_{1}),e_{2}\rangle=d\langle e_{1},e_{2}\rangle.\]
As usual, a $*$-connection is \emph{Hermitian} when $\nabla=\nabla^{*}$, i.e.
\[\langle e_{1},\nabla(e_{2})\rangle - \langle \nabla(e_{1}),e_{2}\rangle=d\langle e_{1},e_{2}\rangle.\]
\newline\newline
We call two $C^{k}$-modules $\mathpzc{E},\mathpzc{F}$ \emph{topologically isomorphic} if there exists an invertible adjointable operator $g:E^{k}\rightarrow F^{k}$. Such $g$ extends to a topological isomorphism between $E^{i}$ and $F^{i}$ for all $i\leq k$.
\begin{theorem}[\cite{Mes}] Let $B$ be a $C^{k}$-algebra, $\mathpzc{E}\lrh B$ a $C^{k}$-module, and $(\mathpzc{F},T)$ a transverse $C^{k}$ $KK$-cycle for $(B,C)$. If $\nabla:E^{k}\rightarrow E^{k}\tildeotimes_{\mathcal{B}_{k}}\Omega^{1}(\mathcal{B}_{k})$ is a Hermitian connection, then the operator
\[1\otimes_{\nabla} T:E^{k-1}\otimes\Dom T\rightarrow E^{k-1}\tildeotimes_{\mathcal{B}_{k-1}}F^{k-1},\]
is essentially selfadjoint and regular in $E^{k-1}\tildeotimes_{\mathcal{B}_{k-1}}F^{k-1}$. Morever, the graphs \[\mathfrak{G}((1\otimes_{\nabla}T)_{i})^{j}\subset E^{j}\] are topologically isomorphic to $E^{k}\tildeotimes_{\mathcal{B}_{k}}\mathfrak{G}(T_{i})^{j}$, for $i+j\leq k$. 
\end{theorem}
The operator $1\otimes_{\nabla} T$ is symmetric because $\nabla$ is Hermitian. Note that in this theorem the transversality property enters to make sure that each $\mathfrak{G}(T_{i})^{j}$ is a left $\mathcal{B}_{k}$-module for $i+j\leq k$. Also, it should be noted that the isomorphism $\mathfrak{G}(1\otimes_{\nabla}T)_{i})^{j}\rightarrow E^{i}\tildeotimes\mathfrak{G}(T_{i})$ is the identity in the first coordinate. As such it gives a description of the domain of the operator $(1\otimes_{\nabla}T)_{i}^{j}$, see \cite{Mes} for details. Transverse smoothness of connections is defined straightforwardly, again in an inductive way.
\begin{definition} A connection $\nabla:E^{k}\rightarrow E^{k}\tildeotimes_{\mathcal{B}_{k}}\Omega^{1}(\mathcal{B}_{k})$ on a $C^{k}$-cycle $(\mathpzc{E},D)$ is said to be a transverse $C^{k}$-connection if $[D,\theta_{i}(\nabla)]$ extends to a completely bounded operator $\mathfrak{G}(D_{i})^{j}\rightarrow\mathfrak{G}(D_{i})^{j}\tildeotimes_{\mathcal{B}_{j}}\Omega^{1}(\mathcal{B}_{j})$ for all $i+j\leq k$. Equivalently, if
\[(\textnormal{ad}(D))^{n}(\nabla)(D\pm i)^{-n+1},\quad (D\pm i)^{-n+1}(\textnormal{ad}(D))^{n}(\nabla),\]
extend to completely bounded operators $E^{j}\rightarrow E^{j}\tildeotimes_{\mathcal{B}_{j}}\Omega^{1}(\mathcal{B}_{j})$, for $n+j\leq k$.
\end{definition}
Note that a transverse $C^{k}$ connection induces connections $\theta_{i}(\nabla):\mathfrak{G}(D_{i})^{j}\rightarrow \mathfrak{G}(D_{i})^{j}\tildeotimes_{\mathcal{B}_{j}}\Omega^{1}(\mathcal{B}_{j})$ for all $i+j\leq k$. These connections are not Hermitian for the inner product on $\mathfrak{G}(D_{i})^{j}$, but they are $*$-connections.

\begin{definition} Let $A,B$ be $C^{k}$-algebras. A \emph{geometric correspondence} is a $C^{k}$-cycle with connection. That is, it is a triple $(\mathpzc{E},D,\nabla)$, where $\mathpzc{E}$ is a $C^{k}$-bimodule, $D$ a $C^{k-1}$ operator, and $\nabla$ a transverse $C^{k}$-connection.
\end{definition}

\section{The product construction}
Geometric correspondences can be composed according to the algebraic procedure described in Proposition \ref{alg}. The smoothness conditions imposed on the cycles make sure that this algebraic procedure preserves all the desired analytic properties. In particular, the smoothness conditions themselves are preserved.
\begin{theorem}[\cite{Mes}]\label{product} Let $A,B,C$ be $C^{k}$-algebras, with $k\geq 2$, and $(\mathpzc{E},S,\nabla)$, $(\mathpzc{F},T,\nabla ')$ $C^{k}$-cycles with connection. Then
\[(\mathpzc{E}\tildeotimes_{B}\mathpzc{F},S\otimes 1+ 1\otimes_{\nabla} T, 1\otimes_{\nabla}\nabla '),\]
is a $C^{k}$-cycle with connection. It represents the Kasparov product of $(\mathpzc{E},S)$ and $(\mathpzc{F},T)$.
\end{theorem}
\begin{remark} The condition $k\geq 2$ is needed to guarantee that the connection on the module is again trnasverse $C^{k}$. If one just wants to compute Kasaprov products, one can work with $C^{1}$-modules. 
Commutator conditions are direct. A result of Kucerovsky \cite{Kuc} on unbounded Kasparov products then gives the last assertion.
\end{remark}
We can view geometric correspondences as morphisms of spectral triples. A morphism between $C^{k}$ spectral triples $(A,\mathpzc{H},D)$ and $(B,\mathpzc{H}',T)$ is a $C^{k}$-bimodule with connection $(\mathpzc{E},S,\nabla)$ such that the spectral triple
\[(A,\mathpzc{E}\tildeotimes_{B}\mathpzc{H}', S\otimes 1 +1\otimes_{\nabla}T),\]
is $C^{k}$ unitarily isomorphic to $(A,\mathpzc{H},D)$. There is a category of spectral triples for each degree of smoothness. If we denote the category of $k$-smooth spectral triples by $\Psi^{k}$, then Theorem \ref{product} says that the bounded transform
\[\mathfrak{b}:(\mathpzc{E},D,\nabla)\mapsto [(\mathpzc{E},D(1+D^{2})^{-\frac{1}{2}})],\]
is a functor $\Psi^{k}\rightarrow \mathfrak{KK}$.

\end{document}